\documentclass[3p,12pt]{elsarticle}
\usepackage{amsthm}
\usepackage{amssymb}
\usepackage{amsmath}
\usepackage{epstopdf}
\usepackage{graphicx}
\usepackage{multirow}
\usepackage{flafter}
\usepackage[section]{placeins}
\usepackage{booktabs}
\usepackage{array}
\usepackage{tabularx}
\usepackage[ruled]{algorithm2e}
\usepackage{color}

\renewcommand{\arraystretch}{1.12}
\numberwithin{equation}{section}

\newproof{pf}{Proof}[section]
\newdefinition{rmk}{Remark}

\biboptions{sort&compress}

\journal{}

\begin{document}

\begin{frontmatter}

\title{A Two-Stage Optimization Framework for Real-World Large-Scale University Course Timetabling Based on Orthogonalized Hierarchical Modeling and Conflict-Driven ALNS Repair}

\author[A]{Han-Bin Liu}
\author[A,B,C]{Huayi Wei}
\author[A]{Yanlin Li\corref{cor1}}
\ead{liyanlin@xtu.edu.cn}
\cortext[cor1]{Corresponding author.}

\address[A]{School of Mathematics and Computational Science, Xiangtan University, Xiangtan 411105, China}
\address[B]{National Center of Applied Mathematics in Hunan, Xiangtan University, Xiangtan 411105, China}
\address[C]{Hunan Key Laboratory for Computation and Simulation in Science and Engineering, Xiangtan University, Xiangtan 411105, China}

\begin{abstract}
University course timetabling is a large-scale combinatorial optimization problem involving strongly coupled temporal and resource constraints.
This paper proposes a two-stage optimization framework for real-world large-scale university course timetabling, integrating orthogonalized hierarchical modeling with conflict-driven adaptive large neighborhood search (ALNS) repair.
The first stage performs global course-level optimization, while the second identifies residual conflicts through a conflict graph and localizes them into connected conflict regions.
The Course--Occurrence--Session representation provides the structural basis for multi-granularity repair, which employs Cluster-, Course-, Occurrence-, and Session-level destroy operators within the ALNS framework, followed by fine-grained greedy repair.
The framework is evaluated on a real-world university timetabling instance involving 3,073 course groups, 8,015 scheduling events, 998 student classes, and 1,511 instructors.
Ablation experiments demonstrate the effectiveness of hierarchical repair granularity and adaptive operator selection.
Across ten independent runs, the ALNS + Greedy framework achieves a mean residual conflict count of 18.50, compared with 172 for deterministic Greedy-only repair under the same initial scheduling conditions, representing an 89.24\% reduction.
These results support the effectiveness of the proposed conflict-driven hierarchical repair framework for large-scale university course timetabling.
\end{abstract}

\begin{keyword}
Large-Scale University Course Timetabling \sep Metaheuristic Optimization \sep Hierarchical Modeling \sep ALNS \sep Hierarchical Conflict Repair
\end{keyword}

\end{frontmatter}

\section{Introduction}

The University Course Timetabling Problem (UCTP) is a challenging combinatorial optimization problem that requires assigning courses to feasible time periods and teaching resources while satisfying a diverse set of institutional requirements \cite{Schaerf1999,BurkePetrovic2002,Lewis2008}. In practical university environments, scheduling decisions must coordinate courses, instructors, student classes, rooms, and temporal requirements over a multi-week horizon. Additional complications arise from institution-specific teaching arrangements, resource-capacity restrictions, and pre-assigned courses whose schedules must remain unchanged. These requirements create strong dependencies among scheduling decisions: modifying one assignment may affect multiple instructors, student classes, and related courses. As the number of scheduling entities and constraints increases, constructing a high-quality timetable while controlling constraint violations becomes computationally demanding \cite{KristiansenStidsen2013,Chen2021}.

A broad range of methods has been investigated to address the computational complexity of UCTP, including mathematical programming, local search, population-based metaheuristics, hybrid optimization, and problem decomposition \cite{BurkePetrovic2002,Lewis2008,KristiansenStidsen2013}. Decomposition-based methods seek to reduce the effective search space by separating a complex scheduling problem into more manageable subproblems or stages. For example, De Causmaecker et al. \cite{DeCausmaecker2009} developed a decomposed metaheuristic for a real-world university timetabling problem, while Almeida et al. \cite{Almeida2023} investigated a hybrid metaheuristic approach incorporating instance decomposition for large-scale course timetabling. These studies demonstrate the value of combining structured problem organization with heuristic search. Nevertheless, decomposition of the original problem and the subsequent repair of conflicts remaining in an optimized timetable represent different tasks. The former primarily aims to facilitate the initial optimization process, whereas the latter must determine which parts of an existing timetable should be reconsidered without unnecessarily disturbing assignments that are already satisfactory.

Large Neighborhood Search (LNS) and Adaptive Large Neighborhood Search (ALNS) provide a complementary mechanism for exploring complex scheduling configurations through destroy-and-repair operations \cite{Shaw1998,RopkePisinger2006}. By temporarily releasing a subset of scheduling decisions and reconstructing them, these methods can explore neighborhoods that are difficult to reach through elementary local moves. Kiefer et al. \cite{Kiefer2017} demonstrated the effectiveness of ALNS for curriculum-based course timetabling by combining problem-specific destroy-and-repair operators with adaptive operator selection. 
Related research has also investigated adaptive neighborhood selection and competition-guided search strategies to allocate search effort among alternative neighborhoods \cite{SongChenXuWangSongTang2020}.
These studies establish that the design and selection of search operators are important to timetabling performance. However, selecting among operators and defining the structural scope of the assignments affected by those operators are distinct design questions. In particular, when global optimization leaves a structured set of residual conflicts, the repair process must determine not only which operator to apply, but also which scheduling entities should be reconsidered jointly and which assignments should remain invariant.

This issue is closely related to how scheduling entities are represented. Existing research has explored alternative temporal and structural representations of course timetabling. For example, course-pattern formulations represent feasible temporal arrangements through patterns associated with courses \cite{Bagger2019}, while customized multi-level timetabling methods organize scheduling decisions according to different planning structures and criteria \cite{DunkeNickel2023}. Such approaches provide important foundations for representing temporal dependencies and accommodating complex institutional requirements. However, the semantic organization of scheduling entities and the definition of repair neighborhoods serve related but non-identical purposes. A representation suitable for describing a complete course or its temporal patterns does not, by itself, specify whether a conflict should be repaired by rescheduling an entire course, modifying one of its weekly occurrences, or changing only a temporally localized part of its teaching arrangement. For practical multi-week timetabling, this distinction is important because a course may contain several recurring weekly occurrences, while the teaching attributes associated with those occurrences may change over different intervals of the semester.

Motivated by these considerations, this study develops an orthogonalized hierarchical modeling approach that distinguishes three complementary semantic responsibilities: Course represents the complete teaching entity, Occurrence preserves the identity of a recurring within-week teaching event, and Session represents a maximal contiguous interval over which the relevant teaching attributes remain stable. These entities are formally related, but they are not interchangeable: the identity of a weekly occurrence is preserved across changes in teaching attributes, whereas a session captures the temporal interval over which those attributes remain unchanged. The term orthogonalized refers to this separation of modeling responsibilities rather than mathematical orthogonality or complete independence among the entities. The resulting representation provides a structured basis for identifying the scheduling units that may need to be reconsidered at different levels of granularity.

The proposed modeling approach is further connected to a conflict-driven repair mechanism. Following global metaheuristic optimization, residual violations are identified from the decoded timetable and represented through a conflict graph. The connected components of this graph define localized conflict regions, referred to as Clusters, which determine where subsequent repair effort is concentrated. Cluster is therefore a conflict-derived repair region rather than an additional semantic level in the Course--Occurrence--Session representation. Within these regions, four destroy operators release adjustable scheduling decisions at the Cluster, Course, Occurrence, and Session levels, respectively. Explicit repair-scope constraints preserve assignments outside the selected scope, while the repair procedure reconstructs the released decisions subject to the applicable scheduling constraints. This design establishes a direct connection between the semantic organization of scheduling entities, the localization of residual conflicts, and the granularity of the search neighborhood.

The methodological contribution of this study is therefore not simply the combination of a metaheuristic optimizer and ALNS, nor the independent introduction of multiple destroy operators. Instead, it lies in organizing the repair process around the interaction between an orthogonalized hierarchical representation of scheduling entities and the structure of residual conflicts in the current timetable. The resulting two-stage framework first performs global course-level optimization and then applies conflict-driven hierarchical repair, followed by fine-grained greedy repair of the remaining violations. This coarse-to-fine strategy is intended to provide sufficient freedom to reorganize strongly coupled assignments while limiting unnecessary modifications to parts of the timetable outside the selected repair scope.

The proposed framework is evaluated on a large-scale real-world university timetabling instance comprising 3,073 course groups, 8,015 scheduling events, 998 student classes, and 1,511 instructors. The instance involves heterogeneous temporal and instructional structures, making it suitable for assessing a framework that must coordinate global search with structured local repair. Comparative experiments and ablation studies are conducted to investigate the effectiveness of the overall repair process, the contribution of different repair granularities, and the influence of adaptive operator selection.

The main contributions of this work are summarized as follows:
\begin{itemize}
    \item \textbf{Orthogonalized hierarchical modeling.}
    A formal Course--Occurrence--Session representation is developed to distinguish the complete course entity, recurring within-week occurrence identity, and cross-week intervals of stable teaching attributes. This representation separates their semantic responsibilities while retaining the mappings needed for timetable evaluation and repair.

    \item \textbf{Conflict-driven localization and multi-granularity repair.}
    A residual conflict graph is constructed after global optimization, and its connected components are used to identify localized repair regions. Four destroy operators operate at the Cluster, Course, Occurrence, and Session levels, with explicit scope constraints preserving assignments outside the selected repair scope.

    \item \textbf{A two-stage optimization and hierarchical repair framework.}
    Global metaheuristic optimization is integrated with adaptive large-neighborhood repair and subsequent fine-grained greedy repair. The framework connects global exploration, residual-conflict localization, and multi-granularity neighborhood search within a unified coarse-to-fine procedure.

    \item \textbf{Evaluation on a large-scale real-world instance.}
    The proposed framework is evaluated on a university timetabling instance involving 3,073 course groups, 8,015 scheduling events, 998 student classes, and 1,511 instructors. Comparative experiments and ablation studies assess the effectiveness of the repair framework and the contributions of its principal components.
\end{itemize}

The remainder of this paper is organized as follows. Section~\ref{Related Work} reviews research on structured modeling, problem decomposition, conflict-oriented search, and large-neighborhood methods for university course timetabling. Section~\ref{Problem Formulation} formulates the considered timetabling problem and defines the scheduling entities, constraints, and objective function. Section~\ref{proposed method} presents the proposed two-stage optimization and hierarchical conflict repair framework. Section~\ref{Experimental Results and Analysis} describes the experimental configuration and reports the comparative and ablation results. Finally, Section~\ref{Conclusion} concludes the paper and discusses potential directions for future research.

\section{Related Work}\label{Related Work}

University course timetabling is a constrained combinatorial optimization problem involving temporal assignments, resource availability, and institution-specific requirements. Existing research has investigated metaheuristic optimization, large-neighborhood search, problem decomposition, and structural modeling. This section reviews these directions to clarify the methodological context of the proposed framework.

\subsection{Metaheuristic Optimization and Problem Decomposition}

Metaheuristic methods are widely used for university course timetabling because of the difficulty of satisfying interacting constraints while maintaining timetable quality. Previous studies have investigated local search, evolutionary algorithms, and hybrid strategies \cite{Schaerf1999,BurkePetrovic2002,Lewis2008,Abdullah2007,JatYang2011,ChenShih2013,Rezaeipanah2021,ThepphakornPongcharoen2020}. For large and customized instances, problem decomposition provides an additional means of managing computational complexity. De Causmaecker et al. \cite{DeCausmaecker2009} developed a decomposed metaheuristic for a real-world timetabling problem with overlapping time slots and irregular weekly schedules. Almeida et al. \cite{Almeida2023} combined ALNS, guided local search, and variable-neighborhood search with an instance-decomposition strategy for large-scale course timetabling.

These studies establish the value of combining global search with structured problem reduction. However, instance decomposition and residual-conflict repair address different primary tasks. The former organizes the original problem into subproblems, whereas the latter starts from an existing timetable and identifies assignments that require further modification. The proposed framework focuses on the latter setting, using residual conflicts to guide localized repair after global optimization.

\subsection{Large-Neighborhood Search and Adaptive Search}

Large Neighborhood Search (LNS) improves a current solution by removing selected assignments and reconstructing the resulting partial solution \cite{Shaw1998}. Adaptive Large Neighborhood Search (ALNS) extends this approach by selecting among multiple destroy and repair operators according to their observed performance \cite{RopkePisinger2006}. In course timetabling, Kiefer et al. \cite{Kiefer2017} developed problem-specific ALNS operators with adaptive selection, while Song et al. \cite{SongChenXuWangSongTang2020} investigated competition-guided selection among multiple neighborhood structures. Related matheuristic work has also combined heuristic search with mathematical optimization for customized multi-level, multi-criteria timetabling \cite{DunkeNickel2023}.

These methods demonstrate the importance of neighborhood design and operator coordination. However, operator selection and repair-scope definition are distinct decisions: adaptive selection determines which operator to apply, whereas the repair scope determines which assignments may change. The proposed framework addresses the latter by defining Cluster-, Course-, Occurrence-, and Session-level destroy operators according to the residual conflict structure and the organization of course entities. Assignments outside the selected release scope remain unchanged during reconstruction.

\subsection{Structural Modeling and Conflict-Oriented Repair}

Course timetabling models must account for temporal patterns and heterogeneous teaching requirements. Bagger et al. \cite{Bagger2019} incorporated daily course patterns into the scheduling formulation, while Davison et al. \cite{DavisonKheiriZografos2025} considered hybrid teaching arrangements in a multi-objective model. Such studies illustrate the importance of representing the temporal and instructional characteristics of real-world courses.

Conflict graphs and decomposition techniques provide complementary means of organizing scheduling decisions \cite{DiGaspero2007,DeCausmaecker2009,Almeida2023}. A conflict graph identifies relationships among conflicting entities, but graph connectivity alone does not determine the appropriate granularity of repair. In particular, practical courses may contain multiple recurring weekly occurrences and changes in teaching attributes over the semester. These features motivate distinguishing the overall course entity, recurring weekly occurrences, and maximal contiguous intervals over which teaching attributes remain stable.

The proposed Course--Occurrence--Session representation separates these modeling responsibilities and provides a structural basis for defining different repair scopes. After global optimization, connected components of the residual conflict graph identify localized repair regions, termed Clusters. A Cluster is not an additional level of the course representation; it is a conflict-derived region that may involve multiple courses. The framework therefore combines conflict-based localization with Cluster-, Course-, Occurrence-, and Session-level repair scopes, whose released assignments are explicitly defined.

\subsection{Research Gaps and Positioning of This Study}

The reviewed studies provide established methods for global search, decomposition, adaptive neighborhood selection, and temporal modeling. Two issues remain central to the setting considered here. First, the relationship between course structure and repair granularity needs to be made explicit when a course contains multiple recurring occurrences and changes in teaching attributes over time. Second, identifying residual conflicts does not by itself determine the repair region or the set of assignments that should be released while preserving unaffected parts of a timetable.

To address these issues, this study proposes a two-stage optimization framework that combines global metaheuristic optimization with conflict-driven hierarchical repair. Residual conflicts are localized through connected components of a conflict graph, while an orthogonalized Course--Occurrence--Session representation supports the definition of multiple repair granularities. ALNS coordinates the corresponding destroy operators, followed by fine-grained greedy refinement. The contribution lies in integrating scheduling-entity representation, residual-conflict localization, and explicitly scoped multi-granularity repair, rather than claiming novelty for metaheuristics, ALNS, graph-based conflict analysis, or hierarchical modeling individually.

\section{Problem Formulation}\label{Problem Formulation}

The UCTP considered in this study is formulated as a large-scale constrained optimization problem involving courses, instructors, student classes, room types, and a multi-week scheduling horizon. In practical university timetabling, a course is not necessarily a homogeneous scheduling entity. A single course may contain multiple teaching stages, instructors, and temporal segments, which are further expanded into occurrences and sessions during timetable construction. Therefore, the formulation adopts a course-level decision representation while evaluating the resulting timetable at the session level.

The objective is to construct a timetable that satisfies the required temporal, resource, and scheduling constraints while minimizing the remaining constraint violations and undesirable scheduling preferences. 
The resulting problem is a large-scale constrained optimization problem with strongly coupled temporal and resource requirements. The following sections formalize the scheduling entities, resource relationships, and practical constraints that define the problem to be solved.

\subsection{Problem Description}

Let
$
\mathcal{I}=\{1,2,\ldots,N\}
$
denote the set of courses, where \(N\) is the total number of courses to be scheduled. Let
$
\mathcal{T}=\{1,2,\ldots,N_T\},
\qquad
\mathcal{G}=\{1,2,\ldots,N_G\},
\qquad
\mathcal{R}=\{1,2,\ldots,N_R\}
$
denote the sets of instructors, student classes, and room types, respectively.
The scheduling horizon is represented by a set of teaching weeks
$
\mathcal{W}=\{1,2,\ldots,W\},
$
a set of days
$
\mathcal{D}=\{1,2,\ldots,D\},
$
and a set of daily time slots
$
\mathcal{H}=\{1,2,\ldots,H\}.
$
A scheduling position is therefore determined by a week, a day, and a starting time slot, together with the duration of the corresponding session.

Each course \(i\in\mathcal{I}\) is associated with a set of instructors, a set of student classes, a required room type, an effective teaching-week range, a weekly occurrence structure, and a set of temporal constraints. The course may also have preferred time periods and forbidden periods. In addition, some courses may already have fixed scheduling assignments and must therefore be considered as occupied resources during the scheduling of the remaining courses.

The relationships between courses and their associated classes and instructors are represented by binary incidence matrices. Let
$
\mathbf{A}^{G}\in\{0,1\}^{N\times N_G}
$
denote the course--class incidence matrix and
$
\mathbf{A}^{T}\in\{0,1\}^{N\times N_T}
$
denote the course--instructor incidence matrix, where
$$
A^{G}_{ig}
=
\begin{cases}
1, & \text{if course }i\text{ is attended by class }g,\\
0, & \text{otherwise},
\end{cases}
$$
and
$$
A^{T}_{it}
=
\begin{cases}
1, & \text{if instructor }t\text{ teaches course }i,\\
0, & \text{otherwise}.
\end{cases}
$$

Let \(r_i\in\mathcal{R}\) denote the required room type of course \(i\), and let \(R_r\) denote the available number of rooms of type \(r\). The room requirement is therefore represented at the room-type level rather than by assigning a specific physical room during the global optimization stage.

\subsection{Scheduling Constraints}

A feasible timetable must satisfy a set of temporal, resource, and practical scheduling requirements. The constraints are evaluated on the decoded session-level timetable \(X(\mathbf{x})\).

\subsubsection{Temporal Constraints}

Each session must be assigned within the effective teaching-week range of its associated course:
$$
w_i^{\mathrm{start}}
\leq
w_s
\leq
w_i^{\mathrm{end}},
\qquad
s\in\mathcal{S}(i).
$$
The starting time of a session must also allow its entire duration to fit within the daily scheduling horizon:
$$
h_s+\delta_s-1
\leq H.
$$
Let
$
F^{\mathrm{global}}_{w,d,h}
\in\{0,1\}
$
denote the global forbidden-time indicator, where a value of \(1\) indicates that the corresponding time position cannot be used. A session cannot occupy a forbidden time position:
$$
F^{\mathrm{global}}_{w_s,d_s,h}=0,
\qquad
h\in\mathcal{H}_s.
$$
Course-specific forbidden periods are handled analogously. Let
$
F^{i}_{w,d,h}\in\{0,1\}
$
denote the forbidden-time indicator associated with course \(i\). Then
$$
F^{i}_{w_s,d_s,h}=0,
\qquad
h\in\mathcal{H}_s,
$$
must hold for every session \(s\in\mathcal{S}(i)\).
In addition to forbidden periods, courses may have preferred time positions. These preferences do not necessarily invalidate a timetable but contribute to the soft scheduling penalty defined in the objective function.

\subsubsection{Temporal Regularity Constraints}

In addition to basic temporal feasibility, practical university course timetabling requires the weekly occurrences of a course to be distributed according to reasonable temporal patterns.
These requirements are intended to avoid timetables that satisfy the basic temporal and resource constraints but result in undesirable teaching arrangements.

Let \(q_i\) denote the required number of weekly occurrences of course \(i\).
For courses with three weekly occurrences, the occurrences are required to follow a Monday--Wednesday--Friday pattern.
This requirement reflects the institutional scheduling practice for three-times-per-week courses and ensures a regular distribution of teaching activities throughout the week.

For courses with two weekly occurrences, the two occurrences should be distributed across different teaching days.
In particular, two occurrences of the same course cannot be scheduled on the same day.
When feasible, at least one teaching day should be left between the two occurrences.
Let \(d_{i,1}\) and \(d_{i,2}\) denote the teaching days of the two weekly occurrences of course \(i\). The preferred spacing condition can be expressed as
$$
|d_{i,1}-d_{i,2}|\geq 2,
$$
while the mandatory non-coincidence condition is
$$
d_{i,1}\neq d_{i,2}.
$$
The former represents a preferred temporal regularity, whereas the latter represents a mandatory scheduling requirement.

The distinction between these two requirements is important in large-scale timetable construction.
When sufficient feasible positions are available, the preferred spacing should be maintained to obtain a more regular weekly timetable.
However, when the preferred spacing conflicts with other mandatory constraints, it may be relaxed.
In such cases, two weekly occurrences may be placed on consecutive teaching days, but they are still prohibited from being scheduled on the same day.

Therefore, the scheduling model distinguishes between temporal feasibility and temporal regularity.
Temporal feasibility determines whether a timetable satisfies mandatory time-related constraints, whereas temporal regularity characterizes whether repeated weekly teaching occurrences are distributed in a practically reasonable manner.
This distinction also provides the basis for the hierarchical conflict repair strategy introduced in Section~4, where changes to an individual session must be considered together with the temporal structure of its associated occurrence and course.

\subsubsection{Class and Instructor Constraints}

A student class cannot attend two sessions simultaneously. For any two distinct sessions \(s\) and \(s'\), if they share at least one student class, their occupied temporal intervals must not overlap:
$$
\mathbf{A}^{G}_{ig}
\mathbf{A}^{G}_{i'g}
=
1
\quad\Longrightarrow\quad
\mathcal{H}_s
\cap
\mathcal{H}_{s'}
=
\varnothing
$$
for every common class \(g\), whenever the corresponding sessions are active in the same week.
Similarly, an instructor cannot teach two sessions simultaneously. If sessions \(s\) and \(s'\) are associated with a common instructor, their temporal intervals must not overlap:
$$
\mathbf{A}^{T}_{it}
\mathbf{A}^{T}_{i't}
=
1
\quad\Longrightarrow\quad
\mathcal{H}_s
\cap
\mathcal{H}_{s'}
=
\varnothing
$$
for every common instructor \(t\) under the same active week.
These constraints generate the principal pairwise resource conflicts considered by the conflict detection procedure in Section~4.3.

\subsubsection{Room Capacity Constraints}

Unlike class and instructor conflicts, room constraints are evaluated at the room-type capacity level. Let \(r(s)\) denote the room type required by session \(s\). For a given week \(w\), day \(d\), and time slot \(h\), define
$$
L_{r,w,d,h}
=
\sum_{s\in\mathcal{S}}
\mathbb{I}
\left[
r(s)=r
\land
s\text{ is active at }(w,d,h)
\right],
$$
where \(L_{r,w,d,h}\) is the number of sessions simultaneously requiring room type \(r\).
The room capacity constraint is then
$$
L_{r,w,d,h}
\leq
R_r,
\qquad
r\in\mathcal{R}.
$$
Thus, a room-capacity conflict occurs when the number of simultaneously active sessions requiring a given room type exceeds the available number of rooms of that type. This formulation is consistent with the room-type resource representation used by the scheduling model.

\subsubsection{Pre-assigned and Special Constraints}

Some courses may have pre-assigned scheduling information that must be preserved during the optimization process. Let
$
\mathcal{S}^{\mathrm{fix}}
$
denote the set of sessions associated with pre-assigned schedules. Their occupied time positions are treated as fixed resource usage when evaluating the remaining courses.

For a newly scheduled session \(s\), conflicts with pre-assigned sessions are evaluated with respect to the shared classes, instructors, and room-type capacity. Consequently, the fixed schedule is not simply excluded from the optimization problem; instead, it contributes to the available resource state against which candidate course assignments are evaluated.

The scheduling model may also include special course-specific requirements. For example, courses with particular weekly scheduling patterns or activity-specific rules may impose additional temporal restrictions. Such requirements are incorporated during decoding and constraint evaluation rather than by introducing separate global decision variables.

\subsection{Heterogeneous Teaching Attributes and Resource Requirements}
\label{subsec:heterogeneous_teaching_attributes}

In real-world university course timetabling, a course does not necessarily correspond to a single instructor, a single teaching-hour type, or a uniform classroom requirement.
Some courses involve multiple instructors who are responsible for different teaching contents or stages, while a course may also comprise different types of teaching activities, such as theoretical lectures, listening sessions, and laboratory sessions.
These teaching activities may have different resource requirements: theoretical lectures are generally assigned to regular classrooms, listening sessions may require classrooms equipped with specific facilities, and laboratory sessions need rooms that satisfy the corresponding experimental teaching requirements.

Consequently, the teaching attributes and resource requirements associated with a course may vary across teaching weeks and within-week occurrences.
For example, the instructor or teaching-hour type associated with a weekly occurrence may change over the semester, leading to corresponding changes in the required classroom type.
Such variations introduce additional complexity into timetable construction, since temporal assignments must remain consistent with the teaching attributes and resource requirements of the corresponding activities.

Therefore, a course-level representation alone may be insufficient to explicitly characterize all temporal variations in teaching attributes and their associated resource requirements.
A suitable modeling approach should preserve the global identity of each course while representing its recurring weekly occurrences and the temporal intervals over which teaching attributes remain unchanged.
This structural requirement motivates the orthogonalized hierarchical modeling approach introduced in Section~\ref{subsec:hierarchical-modeling}.

\section{Proposed Method}\label{proposed method}

This section presents a two-stage optimization framework for real-world large-scale university course timetabling, integrating orthogonalized hierarchical modeling with conflict-driven repair. The first stage performs global course-level optimization, while the second identifies and repairs residual conflicts in the resulting timetable. The Course--Occurrence--Session representation provides the structural basis for conflict localization and multi-granularity ALNS repair, followed by fine-grained greedy refinement. The remainder of this section describes the global optimization stage, the hierarchical modeling and repair procedure, and the overall framework.

\subsection{Metaheuristic-Based Global Optimization}\label{subsec:global-optimization}

The first stage searches the global course-level scheduling space. Its purpose is to obtain a timetable with a favorable overall structure and a reduced number of constraint violations, rather than to require the global search to eliminate every conflict. The resulting timetable is subsequently passed to the conflict-driven repair stage.

\subsubsection{Objective Function}\label{subsubsec:objective-function}

Let $\mathbf{x}\in\Omega$ denote a course-level decision vector, where $\Omega$ is the encoded search space. Decoding $\mathbf{x}$ produces a timetable $X(\mathbf{x})$ at the level of concrete scheduling assignments. Let $C_T(\mathbf{x})$, $C_G(\mathbf{x})$, $C_R(\mathbf{x})$, $C_F(\mathbf{x})$, and $C_X(\mathbf{x})$ denote, respectively, the amounts of instructor conflicts, student-class conflicts, room-capacity violations, forbidden-time violations, and conflicts with pre-assigned schedules. The aggregate hard-constraint violation penalty is
\begin{equation}
C(\mathbf{x}) = w_T C_T(\mathbf{x}) + w_G C_G(\mathbf{x}) + w_R C_R(\mathbf{x}) + w_F C_F(\mathbf{x}) + w_X C_X(\mathbf{x}),
\label{eq:hard-constraint-penalty}
\end{equation}
where $w_T,w_G,w_R,w_F$, and $w_X$ are the corresponding penalty weights.

Let $P(\mathbf{x})$ denote the penalty for soft scheduling preferences, such as preferred time positions. The objective function is formulated as
\begin{equation}
\min_{\mathbf{x}\in\Omega} F(\mathbf{x}) = w_{\mathrm{conflict}} C(\mathbf{x}) + P(\mathbf{x}),
\label{eq:objective-function}
\end{equation}
where $w_{\mathrm{conflict}}$ is a large coefficient used to prioritize hard-constraint conflict reduction over soft preferences. Because the objective uses finite penalty weights, this formulation expresses a weighted priority rather than guaranteeing strict lexicographic optimization in all circumstances.

The evaluation pipeline is therefore
\begin{equation}
\mathbf{x}\longrightarrow X(\mathbf{x})\longrightarrow
\big(C_T,C_G,C_R,C_F,C_X,P\big)\longrightarrow F(\mathbf{x}).
\label{eq:fitness-pipeline}
\end{equation}
The optimizer manipulates the course-level vector, whereas constraint evaluation is performed on the decoded timetable. This separation allows the search algorithm to remain independent of the detailed scheduling and conflict-evaluation logic.

\subsubsection{Course-Level Encoding and Decoding}\label{subsubsec:encoding-decoding}

Let $N$ be the number of courses in the scheduling instance. The global optimizer operates on a decision vector $\mathbf{x}=(x_1,\ldots,x_N)$, where each component encodes the scheduling decision associated with one course. A decoding function $\mathcal{D}(\cdot)$ maps the encoded vector to a timetable:
\begin{equation}
X=\mathcal{D}(\mathbf{x}).
\label{eq:decoding}
\end{equation}
For each course, the encoded value is interpreted according to its temporal range and scheduling information. The decoding procedure expands the course-level decision into the corresponding concrete scheduling assignments, retaining the occurrence and session structure required to evaluate temporal and resource constraints.

The input instance contains course attributes, course--class and course--teacher relationships, room-type requirements, temporal restrictions, preference information, and pre-assigned schedules. These data remain fixed during optimization. They are used by the decoder and evaluation procedure rather than being directly manipulated by the metaheuristic optimizer.

\subsubsection{Fitness Evaluation and Global Search}\label{subsubsec:fitness-evaluation}

For each candidate vector $\mathbf{x}$, the evaluation procedure first obtains $X=\mathcal{D}(\mathbf{x})$, then evaluates instructor conflicts, student-class conflicts, room-capacity violations, forbidden-time violations, conflicts with pre-assigned schedules, and soft preferences. The resulting measures are aggregated using Eq.~\eqref{eq:objective-function}, and the fitness value is returned to the optimizer.

The overall computation is summarized as
\begin{equation}
\boxed{\text{Candidate vector}\rightarrow\text{Decoding}\rightarrow\text{Constraint evaluation}\rightarrow\text{Fitness}}.
\end{equation}
The best timetable obtained at the end of global search becomes the input to the second-stage procedure. Residual conflicts are retained as structured information, rather than being represented only by a scalar fitness value, so that their relationships can guide subsequent repair.

\subsection{Orthogonalized Hierarchical Modeling for Course Timetabling}\label{subsec:hierarchical-modeling}

The global optimizer requires a compact course-level decision representation, whereas conflict detection and repair require more detailed information about recurring teaching events and changes in teaching arrangements over the semester. To meet these requirements, the proposed orthogonalized hierarchical model separates three complementary modeling responsibilities: Course represents the complete teaching entity, Occurrence identifies recurring events within the weekly pattern, and Session captures maximal temporal intervals with stable teaching-attribute profiles. These representations describe different aspects of the same course rather than three independent optimization stages. Their separation avoids conflating course identity, recurring-event identity, and temporal segmentation, while preserving the correspondences required for conflict evaluation and multi-granularity repair.

\subsubsection{Course-Level Representation}\label{subsubsec:course-representation}

A course represents a complete teaching entity and is the basic unit of the global optimization stage. Let $\mathcal{I}=\{1,2,\ldots,N\}$ denote the set of courses. A course may involve multiple student classes, instructors, teaching types, weekly occurrences, and temporal requirements. Its course-level decision coordinates the scheduling of the teaching entity as a whole.

\subsubsection{Weekly Occurrence Representation}\label{subsubsec:occurrence-representation}

A weekly occurrence represents one recurring teaching event within a teaching week. For course $c$, let
\begin{equation}
\mathcal{O}(c)=\{o_{c,1},o_{c,2},\ldots,o_{c,q_c}\},
\end{equation}
where $q_c$ is the number of weekly occurrences required by the course. The index $j$ identifies the recurring position of an occurrence within the weekly pattern; it does not identify a particular calendar week, instructor, or teaching type. Thus, for a course meeting twice per week, $o_{c,1}$ refers to its first weekly occurrence across the teaching period and $o_{c,2}$ to its second weekly occurrence.

This persistent identity makes it possible to represent and repair one recurring weekly occurrence without confusing it with a particular week or a temporary instructor assignment. The occurrence representation is therefore distinct from the segmentation of the semester into intervals with stable teaching attributes.

\subsubsection{Session Representation and Temporal Segmentation}\label{subsubsec:session-representation}

A Session is defined by the maximal contiguous interval of teaching weeks over which the teaching-attribute profile of all weekly occurrences of a course remains unchanged. Let $\theta_{c,w,j}$ denote the attributes associated with course $c$, week $w$, and weekly occurrence $j$. These attributes may include the instructor and teaching-hour type, together with other attributes represented in the scheduling data. Define the profile for week $w$ as
\begin{equation}
\boldsymbol{\theta}_{c,w}=\big(\theta_{c,w,1},\theta_{c,w,2},\ldots,\theta_{c,w,q_c}\big).
\end{equation}
A session $s_{c,k}$ corresponds to a maximal contiguous interval $[a_k,b_k]$ for which $\boldsymbol{\theta}_{c,w}$ remains constant. Maximality means that the interval cannot be extended to an adjacent teaching week without changing the profile.

This definition is important for courses whose teaching arrangements change during the semester. A session is not necessarily attached to only one occurrence: it describes a stable interval of the course-wide weekly profile and may therefore cover multiple weekly occurrences. Conversely, one occurrence sequence may span several sessions when the teaching-attribute profile changes over time. Let $\mathcal{S}(c)$ denote the set of sessions obtained by segmenting course $c$ according to this rule.

The course representation can be summarized as
\begin{equation}
\mathcal{M}(c)=\big(\mathcal{O}(c),\mathcal{S}(c),\Phi_c\big),
\label{eq:course-model}
\end{equation}
where $\Phi_c$ records the correspondence between sessions, weekly occurrences, and their teaching attributes over the associated intervals. The Course level describes the complete teaching entity, the Occurrence level identifies recurring weekly events, and the Session level represents maximal stable intervals of the course's teaching-attribute profile.

\begin{figure}[!h]
    \centering
    \includegraphics[width=0.75\textwidth]{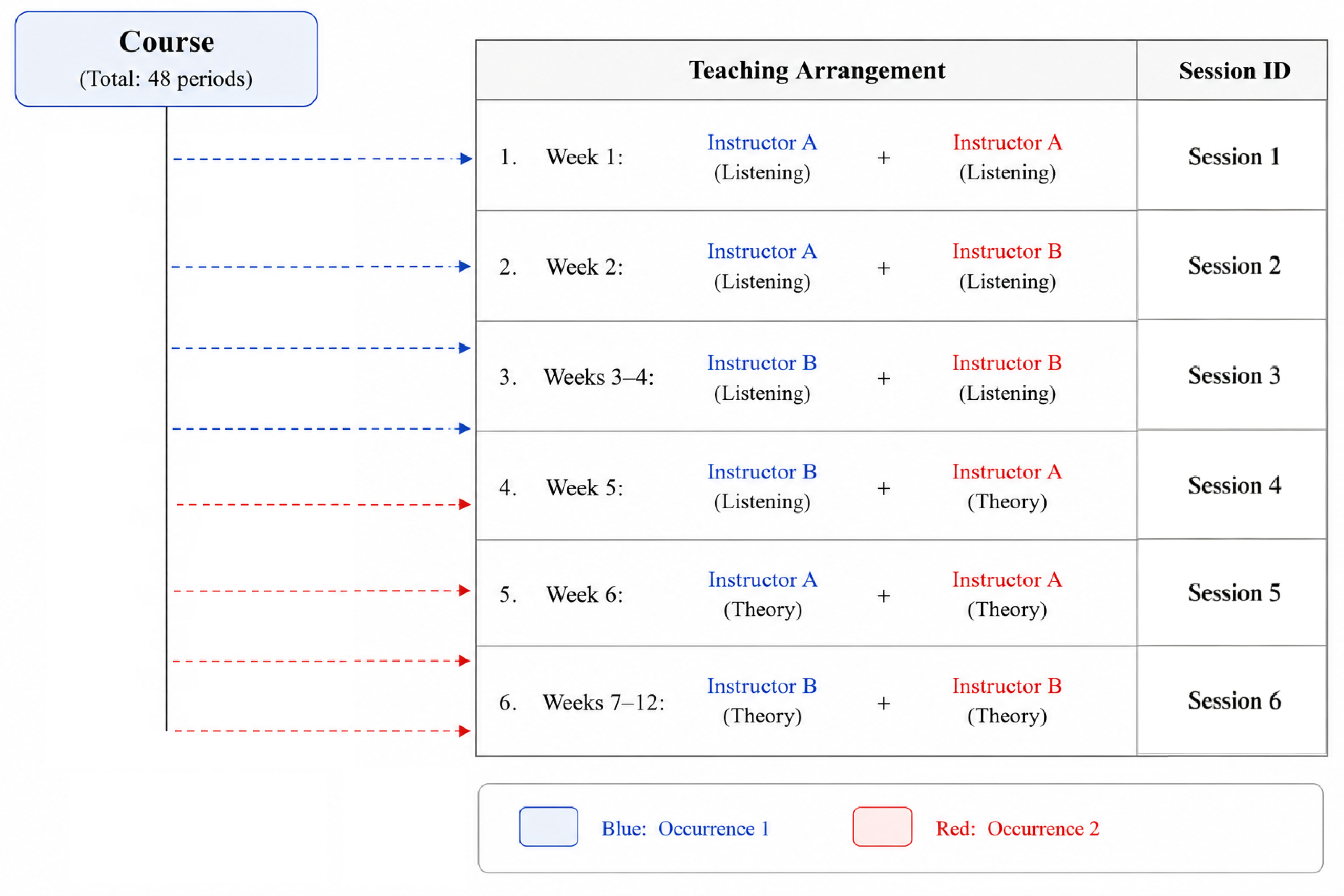}
    \caption{Illustration of the Course--Occurrence--Session representation.}
    \label{fig:course_hierarchy}
\end{figure}

\paragraph{Illustrative example.}
Consider a course with a total workload of 48 periods, consisting of listening and theoretical instruction taught by two instructors. The course meets twice per week, with two periods per meeting, over 12 weeks. The course is represented by one Course entity and two recurring weekly Occurrences. Its semester is then segmented into Sessions according to changes in the joint profile of the two weekly occurrences. For example, a change from one instructor or teaching type to another may start a new session interval; if the profile remains unchanged over several consecutive weeks, those weeks belong to the same session. The resulting segmentation is illustrated in Fig.~\ref{fig:course_hierarchy}.

The representation separates the compact decision unit used for global search from the detailed temporal structure needed for conflict evaluation and repair. In particular, the session segmentation preserves changes in the joint weekly teaching profile that would be lost if a course were treated as temporally homogeneous.

\subsubsection{Relationship Between the Model and Repair Granularity}\label{subsubsec:model-repair-relationship}

The Course--Occurrence--Session representation provides three scheduling-entity granularities, while the conflict graph provides a separate mechanism for identifying a repair region. The cluster is not a fourth level in the course model: it is a connected component of the residual conflict graph and may contain multiple courses. The two structures play complementary roles. A cluster determines where repair is focused; the Course--Occurrence--Session representation determines which part of the scheduling structure can be released within that region.

The framework defines four repair scopes: cluster, course, occurrence, and session. A cluster-level operation may release adjustable assignments associated with courses in a selected conflict component. A course-level operation releases the adjustable schedule of a selected course. An occurrence-level operation releases the assignments associated with a selected recurring weekly occurrence. A session-level operation releases assignments associated with a selected stable interval. These scopes are different operational granularities, not a claim that the corresponding sets are always strictly nested.

For a repair operation $a$, let $\Gamma(a)$ denote its release scope and let $\mathcal{U}(X)$ be the set of basic adjustable scheduling units in timetable $X$. The scope satisfies $\Gamma(a)\subseteq\mathcal{U}(X)$. Assignments outside the selected scope remain unchanged during that operation:
\begin{equation}
x'_u=x_u,\qquad \forall u\in\mathcal{U}(X)\setminus\Gamma(a).
\label{eq:outside-scope-invariance}
\end{equation}
The release scope specifies which assignments may be changed; it does not require every released assignment to receive a different value after reconstruction. Fixed or pre-assigned schedules are excluded from the adjustable release scope.

The representation also supports temporal consistency. Course-level repair can reconsider a course's weekly pattern jointly; occurrence-level repair focuses on a selected recurring event while retaining the other occurrences as references; session-level repair provides a more localized adjustment when the conflict is confined to a stable interval. These operations remain subject to the actual temporal and resource constraints of the instance. Any relaxation of preferred spacing rules must follow the implementation's explicit feasibility policy rather than being assumed as an intrinsic property of the hierarchy.
Thus, the proposed orthogonalization lies in separating the responsibilities of entity representation, recurring-event identification, and temporal segmentation, while linking these representations to the corresponding scopes of global optimization and local repair.

\subsection{Conflict-Driven Multi-Granularity Repair Based on Hierarchical Modeling}\label{subsec:conflict-driven-repair}

After global optimization, the timetable is analyzed to identify residual violations. The repair stage uses structured conflict records to distinguish the types of violations, constructs a graph from pairwise course conflicts, extracts connected components as conflict clusters, and applies ALNS with destroy scopes aligned with the hierarchical scheduling representation. A final fine-grained greedy procedure then targets conflicts that remain after ALNS.

\subsubsection{Residual Conflict Detection and Classification}\label{subsubsec:conflict-detection}

\begin{figure}[!h]
    \centering
    \includegraphics[width=\textwidth]{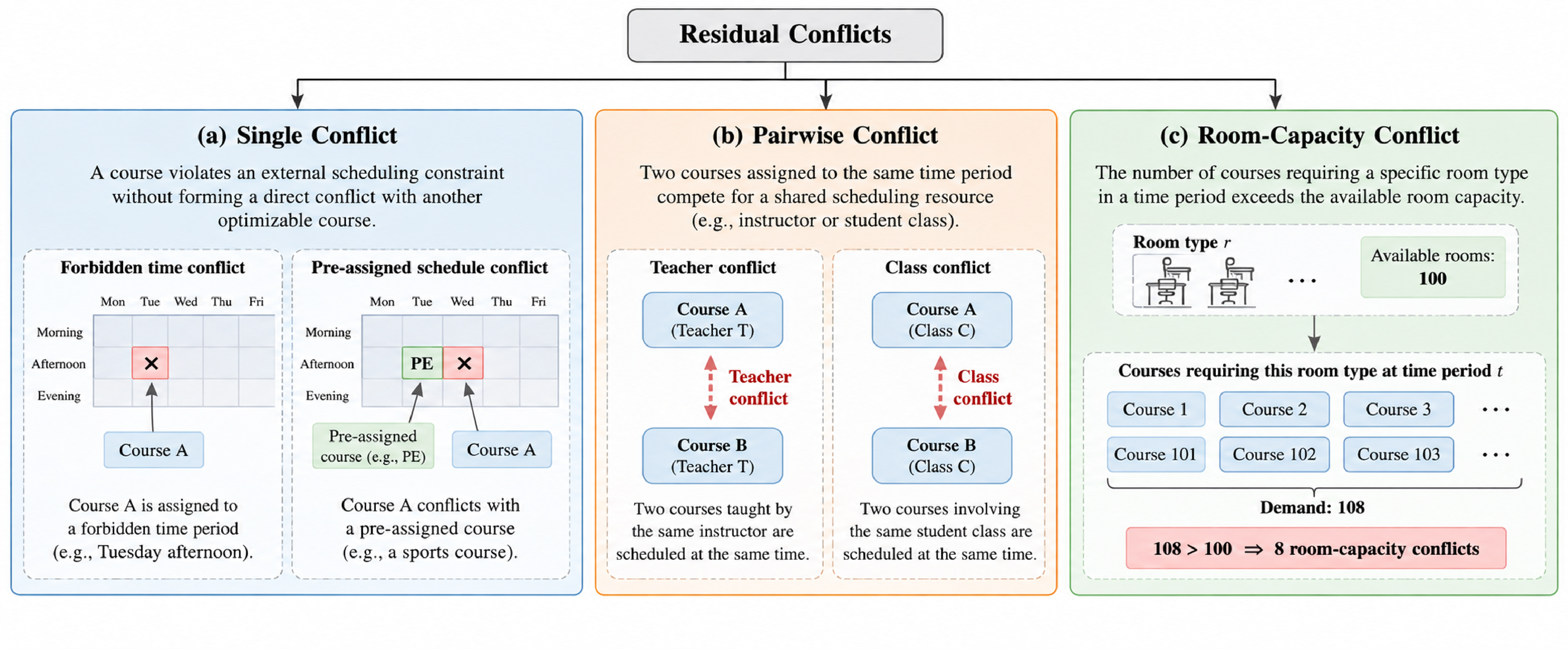}
    \caption{Classification of residual conflicts according to their structural relationships.}
    \label{fig:conflict-classification}
\end{figure}

Residual conflicts are classified according to their structural relationships into single conflicts, pairwise conflicts, and room-capacity conflicts. This classification retains information needed for repair and avoids treating all violations as interchangeable scalar penalties. The categories are illustrated in Fig.~\ref{fig:conflict-classification}.

\paragraph{Single conflicts}
A single conflict occurs when a course violates a constraint without requiring a direct conflict relation with another adjustable course. Examples include an assignment to a forbidden time period or an overlap with a pre-assigned schedule. In the latter case, the pre-assigned schedule is treated as fixed data; the conflict is associated with the adjustable course rather than represented as a conflict between two freely adjustable courses.

\paragraph{Pairwise conflicts}
A pairwise conflict occurs when two courses overlap in time while competing for a shared resource. Typical examples include two courses taught by the same instructor or attended by the same student class at the same time. These conflicts define direct relations between two courses and provide the principal edge information for the residual conflict graph.

\paragraph{Room-capacity conflicts}
A room-capacity conflict is an aggregate overload. Let $D_{r,t}$ denote the number of courses requiring room type $r$ at time $t$, and let $K_{r,t}$ denote the available capacity for that room type and time. The overload is
\begin{equation}
C_R(r,t)=\max\{0,D_{r,t}-K_{r,t}\}.
\label{eq:room-overload}
\end{equation}
For example, if 108 courses require a room type at a time when 100 such rooms are available, the overload is 8. This violation does not necessarily identify a unique pair of courses; it must therefore be retained as an aggregate conflict record rather than automatically reduced to one arbitrary binary edge.

Let $\mathcal{S}^{\mathrm{single}}_{\mathrm{conf}}$, $\mathcal{S}^{\mathrm{pair}}_{\mathrm{conf}}$, and $\mathcal{S}^{\mathrm{room}}_{\mathrm{conf}}$ denote the sessions involved in the respective conflict categories. For a single conflict, the affected session is recorded; for a pairwise conflict, both participating sessions are recorded; and for a room-capacity overload, all sessions contributing to the overload are recorded. The union of these records provides the affected-session information used by subsequent localization and repair.

\subsubsection{Residual Conflict Graph and Cluster Construction}\label{subsubsec:conflict-graph}

Repairing all courses after global optimization would enlarge the neighborhood and could disturb assignments that are not involved in any current violation. The proposed method therefore builds a graph from the actual pairwise conflicts in the optimized timetable. This differs from partitioning a static potential-conflict graph before scheduling: two courses that share a teacher or class need not be connected in the residual graph if their current time assignments do not overlap.

Let $\mathcal{I}^{\mathrm{conf}}$ be the set of courses involved in at least one detected residual conflict. Define the residual conflict graph as $G=(V,E)$, where each vertex represents a course in $\mathcal{I}^{\mathrm{conf}}$ and an edge $(i,j)\in E$ indicates a detected pairwise conflict between courses $i$ and $j$ under the current timetable. Single conflicts and aggregate room-capacity records are retained in the conflict data associated with their affected courses; they should not be interpreted as pairwise edges unless the implementation explicitly constructs such edges.

\begin{figure}[!h]
    \centering
    \includegraphics[width=\linewidth]{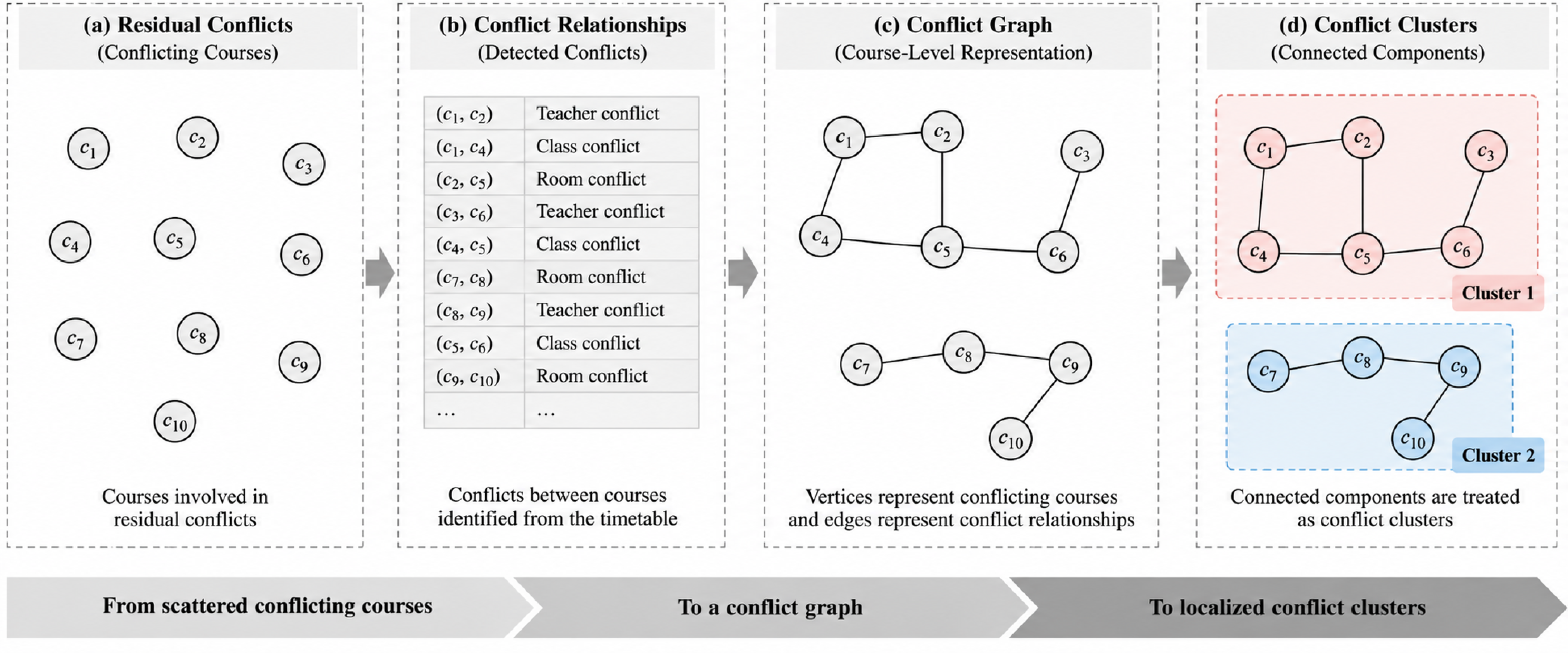}
    \caption{Localization of residual conflicts by graph construction and connected-component decomposition.}
    \label{fig:conflict-clustering}
\end{figure}

The graph is decomposed into connected components $G_1,\ldots,G_K$. Each component defines a conflict cluster $K_k$ with course set $V_k\subseteq V$. Courses in one component are connected through direct or indirect pairwise-conflict paths, whereas no graph edge connects that component to another component. Because the graph is constructed from the current timetable, these clusters reflect active conflict relationships rather than all potential resource relationships in the original instance.

A course with a detected single or aggregate conflict but no pairwise edge must not be discarded. Such a case is retained as a singleton repair cluster or otherwise kept as an explicit repair entry, consistent with the implementation. Thus, the clustering stage preserves all detected conflicts while using connected components to localize pairwise conflict structure.

The resulting transition is
\begin{equation}
\boxed{\text{Decoded timetable}\rightarrow\text{Residual conflicts}\rightarrow\text{Conflict graph}\rightarrow\text{Conflict clusters}}.
\end{equation}
The cluster identifies the region for repair; the hierarchical course model then identifies the possible repair granularity within that region.

\subsubsection{Hierarchical Conflict Representation and Priority Assessment}\label{subsubsec:hierarchical-conflict-priority}

Within each conflict cluster, conflict information is associated with the affected courses, weekly occurrences, and sessions. These representations serve different purposes: course-level aggregation identifies courses with broad involvement in the residual conflict structure; occurrence-level information identifies which recurring weekly events contribute to conflicts; session-level information pinpoints the stable intervals containing affected assignments.

The conflict graph and course hierarchy should therefore be interpreted as complementary structures:
\begin{equation}
\text{Conflict graph}\rightarrow\text{Cluster localization},\qquad
\text{Course}\rightarrow\text{Occurrence}\rightarrow\text{Session}.
\end{equation}
The first determines the repair region, while the second supports increasingly localized analysis and modification within that region.

Conflict degrees are maintained at the cluster, course, occurrence, and session levels to summarize the concentration of residual conflicts. Let $d^{\mathrm{cluster}}$, $d^{\mathrm{course}}$, $d^{\mathrm{occurrence}}$, and $d^{\mathrm{session}}$ denote the corresponding measures. The measures provide priority information for selecting clusters and repair targets. The priority ordering should reflect the conflict records and severity rules implemented by the repair engine; in particular, forbidden-time and fixed-schedule violations may be treated differently from ordinary instructor or class overlaps, while room-capacity violations remain a separate aggregate category.

The resulting priority information guides the choice of repair targets but does not itself change the timetable. Actual changes are made by the destroy and reconstruction operations described below.

\subsubsection{ALNS-Based Multi-Granularity Conflict Repair}\label{subsubsec:alns-repair}

The residual conflicts may differ in both severity and structural coupling. Repairing each conflict independently can cause repeated changes when several conflicts involve related courses. The proposed method therefore uses an adaptive large neighborhood search (ALNS) procedure to reconstruct selected portions of the current timetable.

Let $X$ be the current timetable. For a selected repair scope, a destroy operator releases a subset of adjustable assignments, producing a partial timetable $X^{-}$. A repair operator then reconstructs those assignments and obtains a candidate timetable $X'$:
\begin{equation}
X\xrightarrow{D_k}X^{-}\xrightarrow{R}X'.
\label{eq:destroy-repair}
\end{equation}
The destroy operator determines which assignments may be changed; the repair operator determines how the released assignments are reconstructed. Assignments outside the release scope remain unchanged, as expressed in Eq.~\eqref{eq:outside-scope-invariance}.

The available destroy operators are
\begin{equation}
\mathcal{D}=\{D_{\mathrm{cluster}},D_{\mathrm{course}},D_{\mathrm{occurrence}},D_{\mathrm{session}}\}.
\label{eq:destroy-operator-set}
\end{equation}
They provide different repair granularities. Cluster-based destruction can jointly reconsider multiple courses in a connected conflict region; course-based destruction releases a selected course's adjustable assignments; occurrence-based destruction releases assignments associated with one recurring weekly occurrence; and session-based destruction releases assignments associated with a selected stable interval. The release sets depend on the selected target and its current structure; the four operators should not be assumed to form universally strict nested sets.

\paragraph{Cluster-based destroy}
For a selected cluster $K$ with course set $V_K$, the cluster operator releases the adjustable assignments associated with courses in $V_K$. This allows the repair procedure to reconsider courses connected through chains of active pairwise conflicts within one local neighborhood. Courses outside the selected cluster are not released by this operation.

\paragraph{Course-based destroy}
The course operator selects one or more courses in a conflict cluster and releases their adjustable assignments. Other courses in the cluster remain fixed during that destroy--repair operation. This scope is useful when the cluster contains only a few courses that continue to contribute substantially to the residual conflicts.

\paragraph{Occurrence-based destroy}
For a selected course $c$ and occurrence $o\in\mathcal{O}(c)$, the occurrence operator releases the assignments associated with that recurring weekly event. The other occurrences remain unchanged, providing a more localized neighborhood than releasing the entire course.

\paragraph{Session-based destroy}
The session operator releases assignments associated with a selected session. It provides a fine-grained neighborhood for conflicts concentrated in a limited stable interval. The operation must preserve all structural and feasibility conditions applicable to the released session and leave assignments outside its scope unchanged.

\subsubsection{Greedy Reconstruction of Released Assignments}\label{subsubsec:greedy-reconstruction}

After destruction, the repair operator reconstructs the released units in a sequential, state-dependent manner. Let $\mathcal{U}$ denote the set of units released by the selected destroy operator, and let $\mathcal{P}(u)$ be the candidate scheduling positions for a unit $u\in\mathcal{U}$. The reconstruction procedure orders released units according to their repair priorities, generates candidate positions, filters candidates according to applicable constraints, evaluates the remaining candidates, and assigns a selected position before proceeding to the next unit.

For a candidate position $p$, let $X^{-}_{u\rightarrow p}$ denote the partial timetable after assigning $u$ to $p$. An incremental evaluation can be expressed as
\begin{equation}
\Delta F(u,p)=F(X^{-}_{u\rightarrow p})-F(X^{-}).
\label{eq:incremental-cost}
\end{equation}
When the implementation uses a feasible-candidate set $\mathcal{P}_{\mathrm{feas}}(u)$, the greedy choice is
\begin{equation}
p^*=\arg\min_{p\in\mathcal{P}_{\mathrm{feas}}(u)}\Delta F(u,p).
\label{eq:greedy-position}
\end{equation}
If no candidate satisfies all applicable hard constraints, the implementation's fallback policy determines how the unit is handled; the method description should not imply that a feasible position always exists. After each assignment, resource occupancy and the partial conflict state are updated before the next unit is reconstructed.

The reconstruction respects the granularity of the destroy operation. When an entire course is released, its occurrences and sessions are reconstructed consistently with the course structure. When only an occurrence or session is released, assignments outside that scope are retained. This produces a candidate timetable that can be evaluated by the same conflict-oriented objective used by the global stage.

\subsubsection{Adaptive Operator Selection and Candidate Acceptance}\label{subsubsec:adaptive-selection}

The effectiveness of a destroy operator depends on the current residual-conflict structure. Rather than fixing the operator probabilities throughout the search, the ALNS procedure assigns each operator a weight and selects operators using a roulette-wheel mechanism. Let $w_k^{(t)}$ be the weight of operator $D_k$ at iteration $t$. Its selection probability is
\begin{equation}
p_k^{(t)}=\frac{w_k^{(t)}}{\sum_{j=1}^{m}w_j^{(t)}},\qquad k=1,\ldots,m,
\label{eq:operator-probability}
\end{equation}
where $m=4$ for the four destroy operators in Eq.~\eqref{eq:destroy-operator-set}.

After applying the selected destroy operator and reconstructing its released units, the candidate $X'$ is evaluated. The operator receives a reward according to the search outcome, such as whether the candidate establishes a new best solution, improves the current solution, is accepted without improvement, or is rejected. A generic reward representation is
\begin{equation}
S_k^{(t)}=\begin{cases}
\sigma_1,&\text{if }X'\text{ is a new best solution},\\
\sigma_2,&\text{if }X'\text{ improves the current solution},\\
\sigma_3,&\text{if }X'\text{ is accepted without improvement},\\
\sigma_4,&\text{otherwise},
\end{cases}
\end{equation}
with $\sigma_1>\sigma_2>\sigma_3>\sigma_4\geq0$, provided that these categories and reward values match the implementation.

Let $\bar{S}_k^{(t)}$ denote the average reward of operator $D_k$ during the current adaptation period, and let $\rho\in(0,1]$ be the reaction factor. A typical weight update is
\begin{equation}
w_k^{(t+1)}=(1-\rho)w_k^{(t)}+\rho\bar{S}_k^{(t)},
\label{eq:operator-weight-update}
\end{equation}
followed by recalculation of the probabilities using Eq.~\eqref{eq:operator-probability}. The exact reward accumulation window, initialization, and update timing should be consistent with the implemented algorithm.

Operator selection and candidate acceptance are distinct decisions. An improving candidate can be accepted directly when $F(X')<F(X)$. If the implementation employs a non-greedy acceptance criterion, some non-improving candidates may also be accepted to preserve search diversification. The acceptance rule should be stated in its exact implemented form; the general role of such a rule is to prevent the search from being restricted exclusively to immediate improvements.

The adaptive feedback loop can be summarized as
\begin{equation}
\text{Operator selection}\rightarrow\text{Destroy--repair}\rightarrow\text{Evaluation}\rightarrow\text{Reward and weight update}.
\end{equation}
The adaptive mechanism does not guarantee that broad operators will be selected early and fine-grained operators late. Instead, it adjusts selection probabilities according to observed operator performance, allowing different neighborhood scales to compete throughout the search.

\subsubsection{Fine-Grained Greedy Refinement}\label{subsubsec:fine-greedy-refinement}

After the ALNS phase, the timetable is evaluated again to identify remaining conflicts. A fine-grained greedy refinement procedure then focuses on the affected scheduling units, rather than repeatedly perturbing large conflict regions. For a selected conflicting unit, the procedure generates candidate positions, evaluates their incremental effect on the conflict-oriented objective, applies the best permitted local move, and updates the conflict records. The process continues until the local stopping criterion is reached or no further permitted improvement is found.

This final refinement has a different role from the greedy reconstruction inside each ALNS iteration. During ALNS, greedy reconstruction restores the units released by a selected destroy operator to form a candidate solution. The final refinement, by contrast, is a separate post-ALNS pass intended to eliminate or reduce the conflicts that remain after large-neighborhood search. Both procedures may use incremental evaluation, but their scope and position in the overall algorithm are different.

\subsection{Overall Framework}\label{subsec:overall-framework}

\begin{figure}[!h]
    \centering
    \includegraphics[width=\textwidth]{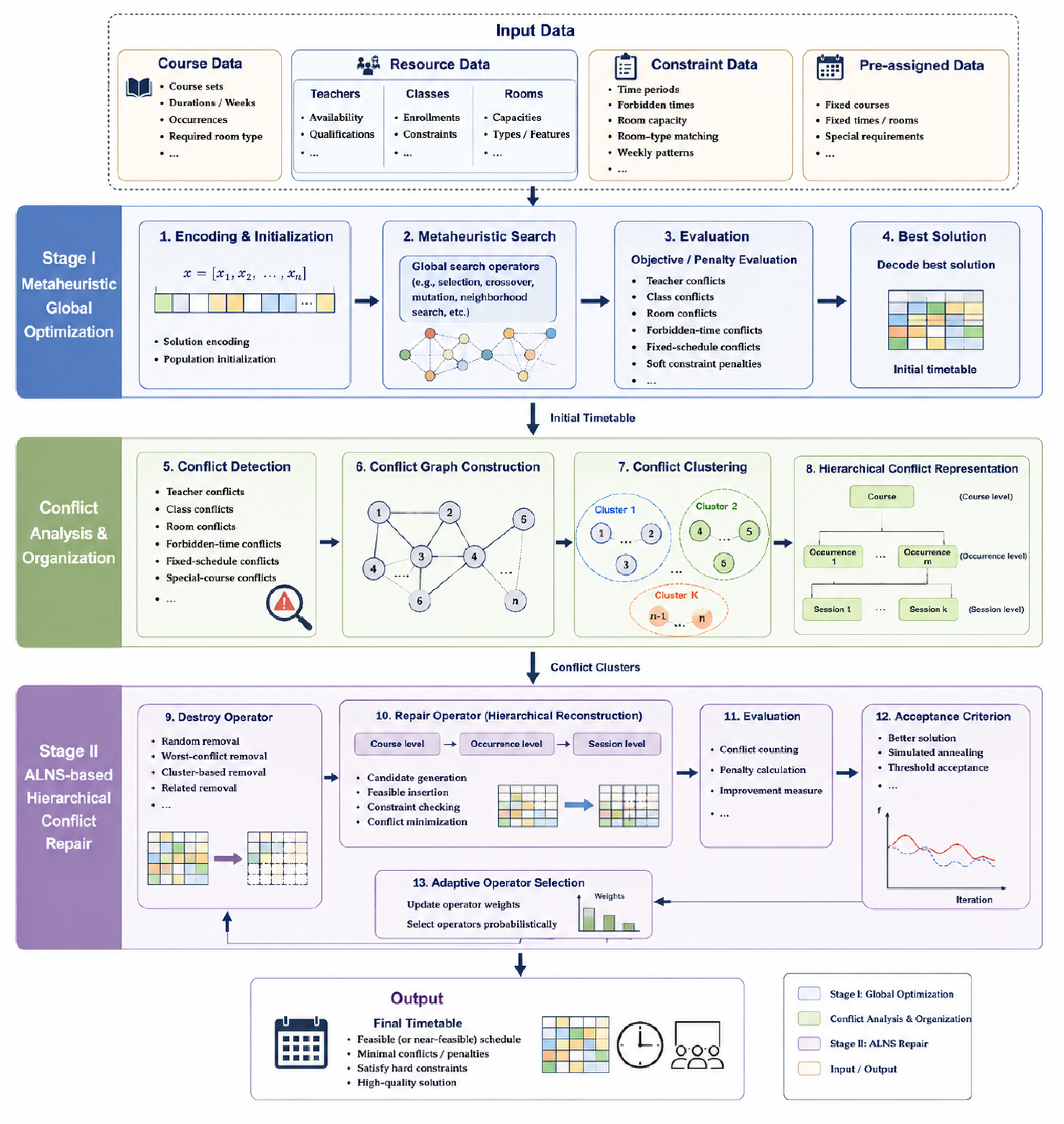}
    \caption{Overall framework of the two-stage optimization and conflict-repair approach.}
    \label{fig:overall-framework}
\end{figure}

The preceding subsections define the global optimization stage, the Course--Occurrence--Session representation, residual conflict localization, and the multi-granularity repair mechanisms. This subsection brings these components together into the complete two-stage framework. The overall workflow is illustrated in Fig.~\ref{fig:overall-framework}.

In Stage 1, the metaheuristic optimizer searches the course-level decision space. Each candidate vector is decoded into a timetable and evaluated using the objective function in Eq.~\eqref{eq:objective-function}. The best timetable found by this stage provides the starting point for residual-conflict analysis.

In Stage 2, the timetable is evaluated to obtain structured records of single, pairwise, and room-capacity conflicts. Pairwise course conflicts are represented in a residual conflict graph, and its connected components define localized conflict clusters. Courses with conflicts that do not produce pairwise edges are retained as singleton clusters or explicit repair entries. The Course--Occurrence--Session representation is then used to identify the appropriate release scope within each repair region. ALNS selects among cluster-, course-, occurrence-, and session-based destroy operators, reconstructs released assignments through greedy repair, evaluates the candidate, and updates operator-selection information according to the adopted adaptive rule. Finally, fine-grained greedy refinement targets the conflicts remaining after ALNS.

The complete process is summarized in Algorithm~\ref{alg:complete-scheduling}. The pseudocode describes the conceptual information flow; the detailed stopping conditions, acceptance rule, reward schedule, and fallback behavior for units with no feasible candidate must match the implementation used in the experiments.

\begin{algorithm}[!h]
\caption{Complete Two-Stage Scheduling Procedure}
\label{alg:complete-scheduling}
\small
\KwIn{Course and resource data, scheduling constraints, pre-assigned schedules, and optimization parameters}
\KwOut{Best timetable $X^*$ and its residual-conflict report}
Initialize the scheduling instance and course-level decision representation\;
Apply the metaheuristic optimizer and decode its best solution to obtain $X$\;
Evaluate $X$ and extract structured residual-conflict records\;
\If{no residual conflicts remain}{
    \Return $X$\;
}
Construct the residual conflict graph from pairwise course conflicts\;
Extract connected components and retain isolated conflict cases as repair entries\;
Build the Course--Occurrence--Session information required for the affected courses\;
Initialize the ALNS destroy-operator weights\;
\While{the ALNS stopping criterion is not satisfied}{
    Select a repair target according to the implemented cluster-selection rule\;
    Select a destroy operator using the current operator probabilities\;
    Release the adjustable assignments within the selected scope to obtain $X^{-}$\;
    Reconstruct the released units using the greedy repair procedure\;
    Evaluate the candidate timetable $X'$\;
    Apply the implemented candidate-acceptance rule\;
    Update the selected operator's performance information and weight\;
    Recalculate the operator-selection probabilities\;
}
Evaluate the current/best timetable and identify remaining conflicts\;
\If{residual conflicts remain}{
    \While{the local stopping criterion is not satisfied}{
        Select a conflicting unit according to the local repair priority\;
        Generate and evaluate candidate positions for the selected unit\;
        Apply the best permitted local move and update conflict records\;
        Stop if no permitted improvement is available\;
    }
}
Return the best timetable $X^*$ and its conflict report\;
\end{algorithm}

In summary, the framework separates global timetable construction from targeted conflict resolution. The global metaheuristic explores the complete course-level search space; residual-conflict analysis identifies the currently coupled repair regions; the hierarchical representation supports multiple repair granularities; ALNS adaptively explores destroy neighborhoods; and final greedy refinement focuses on the remaining localized conflicts. The hierarchy thus functions as the structural link between conflict localization and repair, while the overall algorithm remains a two-stage optimization framework.

\section{Experimental Results and Analysis}\label{Experimental Results and Analysis}

This section evaluates the proposed two-stage optimization framework on a large-scale real-world university course timetabling instance.
The experiments investigate the effects of hierarchical repair granularity and adaptive operator selection, and compare the complete ALNS-based repair framework with the greedy-only baseline.
The dataset characteristics, experimental settings, evaluation metrics, and comparative results are presented and analyzed in the following subsections.

\subsection{Dataset and Problem Characteristics}

The experiments were conducted using a large-scale real-world university course timetabling instance derived from actual academic scheduling data.
Unlike synthetic or benchmark-generated instances, the dataset was obtained from a real university scheduling environment and reflects the structural characteristics and practical constraints encountered in actual course scheduling.
The instance contains 3,073 course groups, which correspond to a 3,073-dimensional course-level decision space in the global optimization stage.
After decoding the course-level decisions according to the Course--Occurrence--Session structure, the resulting timetable contains 8,015 scheduling events that must be evaluated with respect to temporal and resource constraints.
The instance further involves 998 student classes and 1,511 instructors, resulting in substantial coupling among course assignments through shared student, instructor, and temporal resources.
The main characteristics of the real-world experimental instance are summarized in Table~\ref{tab:dataset_characteristics}.

\begin{table}[!h]
\centering
\caption{Characteristics of the large-scale real-world university
timetabling instance derived from actual academic scheduling data.}
\label{tab:dataset_characteristics}
\small
\renewcommand{\arraystretch}{1.15}
\begin{tabular}{lc}
\hline
\textbf{Characteristic} & \textbf{Number} \\
\hline
Course groups & 3,073 \\
Scheduling events & 8,015 \\
Student classes & 998 \\
Instructors & 1,511 \\
\hline
\end{tabular}
\end{table}

\subsection{Experimental Configuration and Protocol}

The global optimization stage employs the Hybrid Crossover Kangaroo Escape Optimizer \cite{LiuWei2026} as the metaheuristic search algorithm.
The population size is set to 300, and the maximum number of iterations is set to 3,000.
The optimizer operates on the course-level decision space, whose dimension is determined by the number of course groups in the experimental instance.
During the global optimization stage, the best-so-far fitness value is recorded at each iteration to monitor the evolution of the scheduling quality.

\begin{figure}[!h]
\centering
\includegraphics[width=0.5\linewidth]{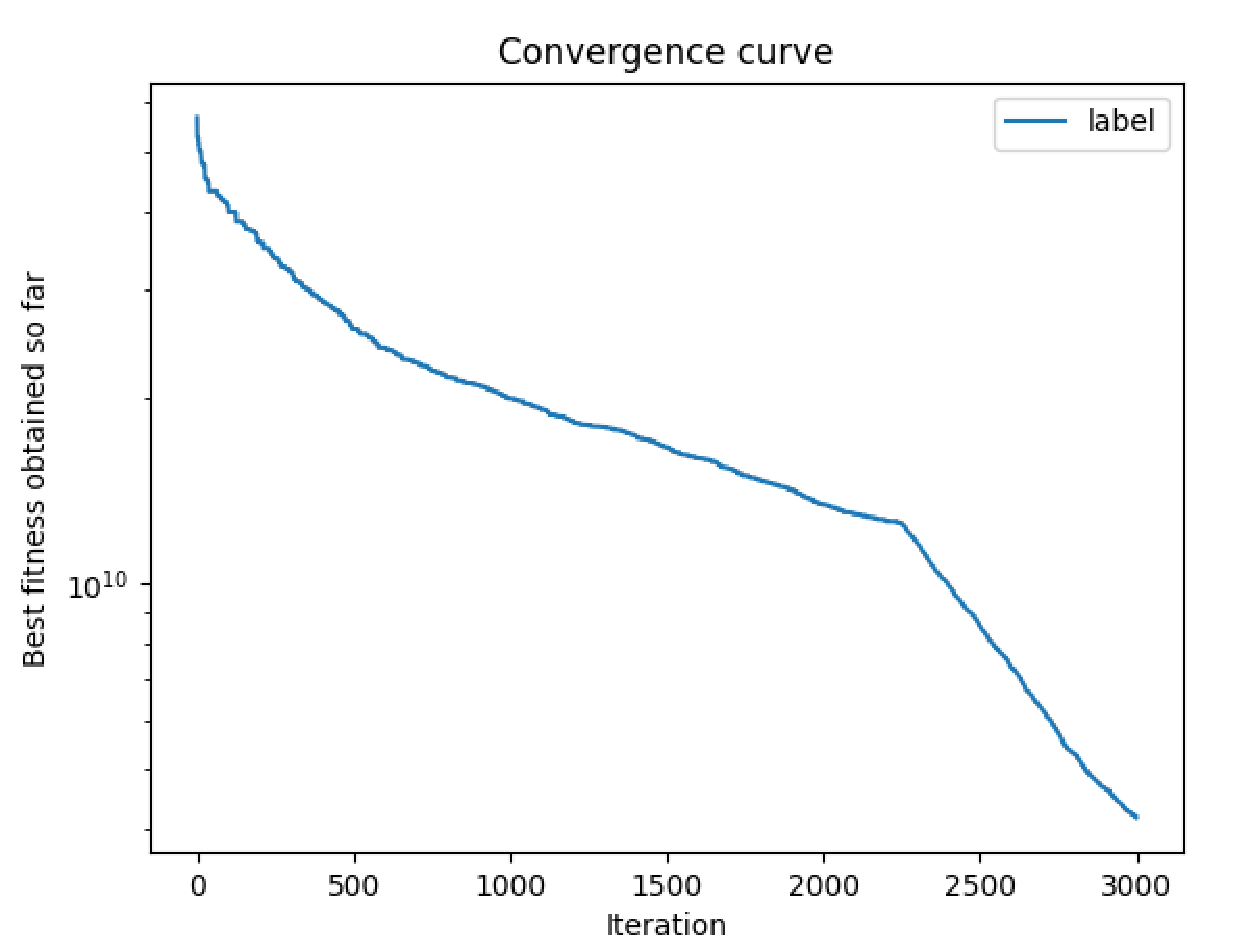}
\caption{Convergence of the global metaheuristic optimization stage over
3,000 iterations. The vertical axis is shown on a logarithmic scale.}
\label{fig:global_convergence}
\end{figure}

After the global optimization stage, the resulting timetable is passed to the proposed conflict detection and hierarchical ALNS-based repair procedure for residual conflict elimination.
The best-so-far fitness value is recorded at each iteration to monitor
the convergence behavior of the global optimization stage.
The resulting convergence curve is shown in
Fig.~\ref{fig:global_convergence}.

The best-so-far fitness decreases throughout the optimization process,
with a relatively rapid improvement during the early iterations
followed by gradual refinement in the subsequent search.
The curve indicates that the global optimization stage continuously
improves the quality of the generated timetable before the solution is
passed to the subsequent conflict repair stage.

\subsection{Evaluation Metrics}

The performance of the proposed scheduling framework is evaluated primarily in terms of the residual conflicts in the resulting timetable.
Since the objective of the scheduling process is to reduce constraint violations, a lower residual conflict count indicates a better scheduling
solution.

Let \(C_{\mathrm{total}}\) denote the total number of residual conflicts.
It is decomposed into the following conflict categories:
\[
C_{\mathrm{total}}
=
C_T+C_G+C_R+C_F+C_X,
\]
where \(C_T\) denotes teacher conflicts, \(C_G\) denotes student-class conflicts, \(C_R\) denotes room-capacity conflicts, \(C_F\) denotes global forbidden-time conflicts, and \(C_X\) denotes conflicts with pre-assigned schedules.
These conflict categories correspond to the main hard constraints considered in the scheduling model and provide a detailed description of the remaining infeasibility after optimization and conflict repair.

The total residual conflict count is used as the primary metric in the ablation experiments.
For each experimental configuration, the final value of \(C_{\mathrm{total}}\) obtained after the specified optimization and repair procedure is recorded.
A smaller value of \(C_{\mathrm{total}}\) indicates that more scheduling conflicts have been successfully eliminated.

For the experiments involving multiple independent runs, the mean, standard deviation, median, minimum, and maximum of the final residual conflict count are reported.
The mean reflects the average performance of the method, while the standard deviation characterizes the stability of the stochastic optimization process.
The median provides a robust measure of the central tendency, whereas the minimum and maximum values represent the best and worst observed outcomes, respectively.

For comparisons between two stochastic optimization configurations, the Wilcoxon signed-rank test is employed to assess whether the observed performance difference is statistically significant.
The test is conducted on the paired results obtained from the same independent runs, with a significance level of \(0.05\).

\subsection{Experimental Configurations and Comparison Settings}

To systematically evaluate the proposed conflict repair framework,
three groups of experiments are conducted.
The first experiment investigates the contribution of the hierarchical
repair granularity by comparing different levels of the
Cluster--Course--Occurrence--Session structure.
The second experiment evaluates the effect of adaptive operator
selection by comparing equal-probability and performance-based adaptive
selection strategies.
The third experiment compares the complete ALNS-based repair framework
with the greedy-only repair strategy used in the previous version of
the scheduling system.

For all experiments, the same real-world university timetabling
instance and the same global optimization procedure are used.
The global optimization stage employs the Hybrid Crossover Kangaroo
Escape Optimizer with a population size of 300 and a maximum of 3,000
iterations.
Unless otherwise specified, the constraint definitions, conflict
detection procedure, evaluation metrics, and common repair mechanisms
are kept consistent across the corresponding experimental
configurations.

\subsubsection{Ablation of Hierarchical Repair Granularity}

The first experiment investigates the contribution of the proposed
hierarchical repair structure.
Three configurations are considered: the complete hierarchical
structure, a configuration without the cluster level, and a
configuration without both the cluster and course levels.

In the Full configuration, the repair process can operate at all four
levels of the proposed hierarchy.
In the w/o Cluster configuration, cluster-level repair is removed while
course-, occurrence-, and session-level repair are retained.
In the w/o Cluster \& Course configuration, both cluster- and
course-level repair are removed, leaving occurrence- and session-level
repair.
All other components of the ALNS framework are kept unchanged so that
the comparison isolates the effect of repair granularity.
As shown in Table~\ref{tab:hierarchical_granularity}, the proposed hierarchical repair structure is progressively simplified by removing the conflict-cluster localization and course-level repair stages.

\begin{table}[!h]
\centering
\caption{Repair structures considered in the ablation study of
hierarchical repair granularity.}
\label{tab:hierarchical_granularity}
\small
\renewcommand{\arraystretch}{1.2}
\begin{tabular}{lc}
\hline
Configuration & Repair Structure \\
\hline
Full
& Cluster $\rightarrow$ Course $\rightarrow$ Occurrence $\rightarrow$ Session \\

w/o Cluster
& Course $\rightarrow$ Occurrence $\rightarrow$ Session \\

w/o Cluster \& Course
& Occurrence $\rightarrow$ Session \\
\hline
\end{tabular}
\end{table}

\subsubsection{Ablation of Adaptive Operator Selection}

The second experiment evaluates the effect of the adaptive operator selection mechanism.
Two configurations are compared: Equal Probability and Adaptive Probability.
Their operator selection strategies are summarized in Table~\ref{tab:operator_selection_config}.

\begin{table}[!h]
\centering
\caption{Configurations for the adaptive operator selection study.}
\label{tab:operator_selection_config}
\small
\renewcommand{\arraystretch}{1.2}
\begin{tabular}{lc}
\hline
\textbf{Configuration} & \textbf{Operator Selection Strategy} \\
\hline
Equal Probability
& Equal selection probability \\
Adaptive Probability
& Performance-based adaptive selection \\
\hline
\end{tabular}
\end{table}

For each configuration, ten independent ALNS repair runs are performed, with 200 ALNS iterations per run.
The same four Destroy operators are used in both configurations, while only their selection strategy is changed.
In the Equal Probability configuration, the four operators are selected with equal probability.
In the Adaptive Probability configuration, the selection probabilities are dynamically updated according to the performance of the operators during the search.
All other experimental settings remain unchanged between the two configurations.

\subsubsection{Comparison between ALNS-based and Greedy-only Repair}

The third experiment evaluates the effectiveness of the proposed two-stage conflict repair framework.
The baseline configuration, denoted as Greedy-only, applies only the fine-grained greedy repair procedure to the residual conflicts of the initial timetable.
The proposed configuration, denoted as ALNS + Fine Repair, first applies ALNS-based large-scale conflict repair and subsequently uses the same fine-grained greedy repair procedure to eliminate the remaining local conflicts.

The comparison is therefore designed to examine whether the introduction of the ALNS-based large-scale repair stage can provide a better starting point for subsequent fine-grained repair.
The same initial timetable and fine-grained repair mechanism are used in both configurations, while the ALNS-based large-scale repair stage is introduced only in the proposed configuration.

\subsection{Results: Hierarchical Repair Granularity}

To investigate the contribution of the proposed hierarchical repair structure, three configurations are evaluated over ten independent runs: Full, w/o Cluster, and w/o Cluster \& Course.
The Full configuration employs the complete Cluster--Course--Occurrence--Session repair hierarchy, whereas the other two configurations progressively remove the coarser repair levels.
The final residual conflict counts are reported in Table~\ref{tab:hierarchical_granularity_results}.

\begin{table}[!h]
\centering
\caption{Run-by-run comparison of different hierarchical repair
granularities over ten independent runs.}
\label{tab:hierarchical_granularity_results}
\small
\renewcommand{\arraystretch}{1.15}
\begin{tabular}{c c c c}
\hline
\textbf{Run}
& \textbf{Full Hierarchy}
& \textbf{Without Cluster}
& \textbf{Without Cluster \& Course} \\
\hline
1  & \textbf{1294} & 1552 & 1692 \\
2  & \textbf{1456} & 1574 & 1536 \\
3  & \textbf{1544} & 1768 & 1610 \\
4  & \textbf{1438} & 1566 & 1736 \\
5  & \textbf{1470} & 1746 & 1534 \\
6  & \textbf{1362} & 1550 & 1514 \\
7  & \textbf{1434} & 1570 & 1594 \\
8  & \textbf{1494} & 1556 & 1670 \\
9  & \textbf{1432} & 1674 & 1594 \\
10 & \textbf{1404} & 1508 & 2088 \\
\hline
Mean
& \textbf{1432.8}
& 1606.4
& 1656.8 \\
Std.
& \textbf{69.26}
& 89.78
& 167.92 \\
Median
& \textbf{1436}
& 1568
& 1602 \\
Min.
& \textbf{1294}
& 1508
& 1514 \\
Max.
& \textbf{1544}
& 1768
& 2088 \\
\hline
\multicolumn{4}{l}{Wilcoxon signed-rank test: Full vs. w/o Cluster,
$p=0.00195$} \\
\multicolumn{4}{l}{Wilcoxon signed-rank test: Full vs. w/o Cluster \& Course,
$p=0.00195$} \\
\hline
\end{tabular}
\end{table}

The Full configuration achieves the lowest average residual conflict count of 1432.8, compared with 1606.4 for w/o Cluster and 1656.8 for w/o Cluster \& Course.
Thus, removing the cluster level increases the mean residual conflict count by 173.6, corresponding to a 10.81\% increase relative to the Full configuration.
Removing both the cluster and course levels further increases the mean residual conflict count to 1656.8, which is 15.63\% higher than that of the Full configuration.

The improvement of the Full configuration is also consistent across the independent runs.
In all ten runs, the Full configuration produces fewer residual conflicts than both w/o Cluster and w/o Cluster \& Course.
The corresponding Wilcoxon signed-rank tests yield \(p=0.00195\) for both comparisons, indicating statistically significant differences at the 0.05 significance level.

In addition to the average performance, the hierarchical structure also improves the stability of the repair process.
The standard deviation increases from 69.26 for the Full configuration to 89.78 when the cluster level is removed, and further to 167.92 when both the cluster and course levels are removed.
The substantially larger maximum value of 2088 in the w/o Cluster \& Course configuration further indicates that restricting the repair process to occurrence- and session-level operations can lead to considerably larger residual conflicts in some runs.

These results demonstrate that the coarse-grained levels of the proposed repair hierarchy provide complementary functionality to the occurrence- and session-level local repair.
The cluster level enables the repair process to jointly reorganize strongly coupled conflict regions, while the course level provides an intermediate repair scope between cluster-level restructuring and fine-grained occurrence/session adjustments.
Removing these coarse-grained levels progressively limits the ability of the repair process to restructure coupled conflicts, resulting in higher residual conflict counts and greater stochastic variation.

To further examine whether the differences observed in the ALNS stage are preserved after the subsequent fine-grained repair, the three hierarchical repair configurations were continued through the same Greedy Repair procedure.
The final residual conflict counts are reported in Table~\ref{tab:hierarchical_repair_final_results}.
The Full Hierarchy configuration achieves a mean of 18.50 residual conflicts over the ten independent runs, compared with 42.80 for the configuration without the Cluster level and 51.10 for the configuration without both the Cluster and Course levels.
The corresponding Wilcoxon signed-rank tests further show significant differences between Full Hierarchy and the two reduced configurations, with $p=0.00391$ and $p=0.00195$, respectively.
Moreover, the Full Hierarchy configuration achieves the lowest residual conflict count in all ten runs, with a maximum of only 34, whereas the two reduced configurations reach maximum values of 66 and 64, respectively.

These results indicate that the hierarchical repair granularity affects not only the intermediate performance of ALNS but also the final outcome after Greedy Repair.
The Cluster and Course levels provide broader neighborhoods for reorganizing structurally related conflicts before the repair process moves toward occurrence- and session-level adjustments.
The consistently lower final conflict counts obtained by the Full Hierarchy configuration suggest that this coarse-to-fine neighborhood structure enables the search to exploit the available repair scope more effectively, thereby providing a more favorable solution for the subsequent Greedy Repair stage.
In contrast, removing the coarse-grained levels restricts the reorganization scope of ALNS, leaving more residual conflicts for the fine-grained repair procedure to resolve.

\begin{table}[!h]
\centering
\caption{Run-by-run comparison of final residual conflict counts
under different hierarchical repair granularities over ten independent runs.}
\label{tab:hierarchical_repair_final_results}
\small
\renewcommand{\arraystretch}{1.15}
\begin{tabular}{c c c c}
\hline
\textbf{Run}
& \textbf{Full Hierarchy}
& \textbf{Without Cluster}
& \textbf{Without Cluster \& Course} \\
\hline
1  & \textbf{8}  & 34 & 44 \\
2  & \textbf{22} & 36 & 64 \\
3  & \textbf{16} & 66 & 62 \\
4  & \textbf{19} & 48 & 56 \\
5  & \textbf{6}  & 60 & 28 \\
6  & \textbf{26} & 31 & 55 \\
7  & \textbf{20} & 34 & 40 \\
8  & \textbf{30} & 40 & 61 \\
9  & \textbf{4}  & 46 & 43 \\
10 & \textbf{34} & 33 & 58 \\
\hline
\textbf{Mean}
& \textbf{18.50}
& 42.80
& 51.10 \\
\textbf{Std.}
& \textbf{10.15}
& 12.09
& 11.75 \\
\textbf{Median}
& \textbf{19.50}
& 38.00
& 55.50 \\
\textbf{Min.}
& \textbf{4}
& 31
& 28 \\
\textbf{Max.}
& \textbf{34}
& 66
& 64 \\
\hline
\multicolumn{4}{l}{Wilcoxon signed-rank test: Full vs. w/o Cluster,
$p=0.00391$} \\
\multicolumn{4}{l}{Wilcoxon signed-rank test: Full vs. w/o Cluster \& Course,
$p=0.00195$} \\
\hline
\end{tabular}
\end{table}

\subsection{Effect of Adaptive Operator Selection}

To investigate the effectiveness of adaptive operator selection in the ALNS-based conflict repair procedure, an ablation experiment was conducted by comparing two operator selection strategies: equal-probability selection and adaptive-probability selection.
Both strategies used the same four destroy operators and were evaluated under the same experimental configuration, with 200 ALNS iterations per run.
Ten independent runs were conducted for each strategy, and the final number of residual conflicts was recorded as the performance measure.
A lower value therefore indicates better conflict repair performance.

The complete run-by-run results are reported in Table~\ref{tab:operator_selection}.
The adaptive-probability strategy achieved a lower residual conflict count in seven out of the ten independent runs.
The equal-probability strategy obtained an average residual conflict count of 1496.2 with a standard deviation of 95.38, whereas the adaptive-probability strategy reduced the average to 1432.8 with a standard deviation of 69.26.
Thus, adaptive operator selection achieved an average reduction of 63.4 residual conflicts, corresponding to a relative improvement of 4.24\%.
Moreover, the standard deviation decreased by approximately 27.39\%, indicating that the adaptive strategy also provides more stable repair performance across independent runs.

\begin{table}[!h]
\centering
\caption{Run-by-run comparison of equal-probability and adaptive-probability
operator selection over ten independent runs.}
\label{tab:operator_selection}
\small
\renewcommand{\arraystretch}{1.15}
\begin{tabular}{c r r r}
\hline
Run & Equal Probability & Adaptive Probability & $\Delta=E-A$ \\
\hline
1  & 1414 & \textbf{1294} & +120 \\
2  & 1558 & \textbf{1456} & +102 \\
3  & \textbf{1482} & 1544 & -62 \\
4  & 1504 & \textbf{1438} & +66 \\
5  & 1736 & \textbf{1470} & +266 \\
6  & 1438 & \textbf{1362} & +76 \\
7  & \textbf{1422} & 1434 & -12 \\
8  & \textbf{1430} & 1494 & -64 \\
9  & 1480 & \textbf{1432} & +48 \\
10 & 1498 & \textbf{1404} & +94 \\
\hline
Mean   & 1496.2 & \textbf{1432.8} & +63.4 \\
Std.   & 95.38  & \textbf{69.26}  & -- \\
Median & 1481   & \textbf{1436}   & -- \\
Min.   & 1414   & \textbf{1294}   & -- \\
Max.   & 1736   & \textbf{1544}   & -- \\
\hline
\multicolumn{3}{l}{Wilcoxon signed-rank test ($p$-value)}
& \textbf{0.0488} \\
\hline
\end{tabular}
\end{table}

To further assess the statistical significance of the observed performance difference, the Wilcoxon signed-rank test was conducted on the paired results of the ten runs.
As reported in Table~\ref{tab:operator_selection}, the resulting two-sided \(p\)-value is 0.0488, which is below the 0.05 significance level.
This result indicates a statistically significant difference between the two operator selection strategies.

Overall, adaptive operator selection achieved a lower mean residual conflict count than equal-probability selection, reducing the average number of residual conflicts from 1496.2 to 1432.8, corresponding to a 4.24\% improvement.
At the same time, the standard deviation decreased from 95.38 to 69.26, representing a reduction of approximately 27.39\%.
These results indicate that adaptive operator selection improves both the average conflict repair performance and the stability of the ALNS search.
The observed improvement supports the use of adaptive operator selection, in which the selection probabilities of destroy operators are dynamically adjusted according to their observed search performance.

\begin{table}[!h]
\centering
\caption{Run-by-run comparison of final residual conflict counts
under equal-probability and adaptive-probability operator selection
over ten independent runs.}
\label{tab:adaptive_operator_selection_final}
\small
\renewcommand{\arraystretch}{1.15}
\begin{tabular}{c r r}
\hline
\textbf{Run} & \textbf{Equal Probability} & \textbf{Adaptive Probability} \\
\hline
1  & 28 & \textbf{8} \\
2  & 27 & \textbf{22} \\
3  & 26 & \textbf{16} \\
4  & 25 & \textbf{19} \\
5  & 38 & \textbf{6} \\
6  & \textbf{22} & 26 \\
7  & 26 & \textbf{20} \\
8  & 39 & \textbf{30} \\
9  & 34 & \textbf{4} \\
10 & \textbf{29} & 34 \\
\hline
\textbf{Mean}
& 29.40
& \textbf{18.50} \\
\textbf{Std.}
& \textbf{5.70}
&10.15 \\
\textbf{Median}
& 27.50
& \textbf{19.50} \\
\textbf{Min.}
& 22
& \textbf{4} \\
\textbf{Max.}
& 39
& \textbf{34} \\
\hline
\multicolumn{3}{l}{Wilcoxon signed-rank test: Equal vs. Adaptive,
$p=0.01367$} \\
\hline
\end{tabular}
\end{table}

To further investigate the effect of adaptive operator selection, the two operator-selection strategies were evaluated under the same hierarchical repair structure and followed by the same Greedy Repair procedure.
The final residual conflict counts over ten independent runs are reported in Table~\ref{tab:adaptive_operator_selection_final}.
The Adaptive Probability configuration achieves a mean of 18.50 residual conflicts, substantially lower than the 29.40 obtained by the Equal Probability configuration.
It also achieves a lower median residual conflict count of 19.50, compared with 27.50 for Equal Probability, while the maximum residual conflict count decreases from 39 to 34.
The paired Wilcoxon signed-rank test gives $p=0.01367$, indicating a statistically significant difference between the two operator-selection strategies at the 0.05 significance level.

The results indicate that adaptive operator selection can improve the final repair performance by allocating search effort according to the observed effectiveness of different Destroy operators.
Instead of assigning the same probability to all operators throughout the search, the adaptive strategy updates their selection probabilities based on their search contributions, allowing more effective neighborhood structures to receive greater opportunities for subsequent exploration.
This mechanism is particularly relevant to the hierarchical repair framework because different conflict structures may require different repair granularities.
The lower final conflict counts obtained by Adaptive Probability suggest that dynamically adjusting the selection of these neighborhoods enables ALNS to exploit the available repair operators more effectively before the subsequent Greedy Repair stage.
It is worth noting that the standard deviation of Adaptive Probability (10.15) is higher than that of Equal Probability (5.70).
Therefore, the present results support an improvement in final repair quality, but do not indicate a reduction in run-to-run variability.

\subsection{Overall Performance and Stability over Ten Independent Runs}

To further evaluate the overall performance and stochastic stability of the complete ALNS + fine-grained greedy repair framework, ten independent runs were conducted on the same large-scale real-world timetabling instance, since the ALNS procedure involves randomized search.
The final residual conflict count \(C_{\mathrm{total}}\) was recorded for each run.

In contrast, the Greedy-only repair strategy used in the previous version of the scheduling system is deterministic.
Under the same initial timetable and scheduling conditions, repeated executions of the Greedy-only procedure produce the same result.
Therefore, its deterministic result of 172 residual conflicts is reported as a fixed reference value in Table~\ref{tab:overall_performance}.
The results are summarized in Table~\ref{tab:overall_performance}.

\begin{table}[!h]
\centering
\caption{Final residual conflict counts of the ALNS + Greedy framework
over ten independent runs and the deterministic Greedy-only baseline.}
\label{tab:overall_performance}
\small
\renewcommand{\arraystretch}{1.15}
\begin{tabular}{c r r}
\hline
\textbf{Run} & \textbf{ALNS + Greedy} & \textbf{Greedy-only} \\
\hline
1  & 8  & \multirow{10}{*}{172} \\
2  & 22 & \\
3  & 16 & \\
4  & 19 & \\
5  & 6  & \\
6  & 26 & \\
7  & 20 & \\
8  & 30 & \\
9  & 4  & \\
10 & 34 & \\
\hline
\textbf{Mean}
& \textbf{18.50}
& \\
\textbf{Std.}
& 10.15
& \\
\textbf{Median}
& 19.50
& \\
\textbf{Min.}
& 4
& \\
\textbf{Max.}
& 34
& \\
\hline
\end{tabular}
\end{table}

Across the ten independent runs, the proposed ALNS + Greedy framework achieved a mean final residual conflict count of 18.50, with a standard deviation of 10.15.
The median value was 19.50, while the minimum and maximum residual conflict counts were 4 and 34, respectively.
In comparison, the deterministic Greedy-only strategy produced 172 residual conflicts on the same instance.
Compared with this reference result, the mean residual conflict count of the proposed framework was reduced by 89.24\%.
Furthermore, the worst result among the ten independent ALNS + Greedy runs was 34 residual conflicts, which remained substantially lower than the 172 conflicts obtained by the Greedy-only strategy.
These results demonstrate the effectiveness of introducing ALNS-based large-scale conflict repair before fine-grained greedy repair.
They also show that the proposed framework maintains a relatively low residual conflict count across the ten independent stochastic runs.

\section{Conclusion}\label{Conclusion}

This paper proposed a two-stage optimization framework for real-world large-scale university course timetabling, integrating orthogonalized hierarchical modeling with conflict-driven adaptive large neighborhood search (ALNS) repair.
The framework combines global timetable optimization with residual-conflict localization through a conflict graph, whose connected components define repair regions.
The Course--Occurrence--Session representation provides the structural basis for four repair granularities---Cluster, Course, Occurrence, and Session---within the ALNS procedure, followed by fine-grained greedy repair.
The key contribution lies in establishing a structured connection between scheduling-entity representation, residual-conflict localization, and multi-granularity repair.

The framework was evaluated on a real-world instance involving 3,073 course groups, 8,015 scheduling events, 998 student classes, and 1,511 instructors.
The hierarchical repair ablation study yielded mean residual conflict counts of 18.50 for the complete configuration, 42.80 without Cluster-level repair, and 51.10 without both Cluster- and Course-level repair; the paired Wilcoxon signed-rank tests indicated statistically significant differences.
Adaptive operator selection achieved a mean of 18.50 residual conflicts, compared with 29.40 for Equal Probability ($p=0.01367$), although its run-to-run variability was higher.
Across ten independent runs, the complete ALNS + Greedy framework achieved a mean of 18.50 residual conflicts, compared with 172 for deterministic Greedy-only repair under the same initial scheduling conditions, corresponding to an 89.24\% reduction.

These results support the effectiveness of integrating hierarchical modeling, conflict-graph localization, and multi-granularity repair within a two-stage optimization framework.
However, the evaluation is based on an instance from a single university.
Future work will investigate the generalizability of the approach across different institutional settings and develop more adaptive repair strategies and efficient incremental and parallel implementations.

\section*{CRediT authorship contribution statement}

\textbf{Han-Bin Liu:} Writing – original draft, Writing – review \& editing, Methodology, Validation. \textbf{Huayi Wei:} Writing – review \& editing, Conceptualization. \textbf{Yanlin Li:} Writing – review \& editing, Methodology.

\section*{Declaration of competing interest}

The authors declare that they have no known competing financial interests or personal relationships that could have appeared to influence the work reported in this paper.

\section*{Ethical and informed consent for data used}
Our research is based on open source data, and there are no human subjects in the article, so ethical and informed consent are not applicable.

\section*{Acknowledgment}

This work was supported by the National Natural Science Foundation of China (Grant Nos. 12371410, 12261131501), the Postgraduate Scientific Research Innovation Project of Hunan Province (Grant No. CX20250944), the Postgraduate Scientific Research Innovation Project of Xiangtan University (Grant No. XDCX2025Y205), and the Construction of Innovative Provinces in Hunan Province (Grant No. 2021GK1010).

\section*{Data availability}

The data that support the findings of this study are available upon request.

\bibliographystyle{elsarticle-num}
\bibliography{reference}

\end{document}